\documentclass[a4paper, 11pt]{amsart}

\usepackage[mathscr]{eucal}
\usepackage{amsmath}
\usepackage{amssymb}
\usepackage{amsfonts}
\usepackage{amsthm}

\usepackage{cases}

\usepackage[mathscr]{eucal}
\usepackage{amsmath}
\usepackage{amssymb}
\usepackage{amsfonts}
\usepackage{amsthm}
\usepackage[dvips]{graphics, color}
\theoremstyle{plain}
\newtheorem{theorem}{Theorem}[section]
\newtheorem{proposition}{Proposition}[section]
\newtheorem{lemma}{Lemma}[section]
\newtheorem{remark}{Remark}[section]

\numberwithin{equation}{section}

\def\<{\left<} \def\>{\right>}

\def\bea{\begin{eqnarray} }
\def\eea{\end{eqnarray} }
\def\be{\begin{equation} }
\def\ee{\end{equation} }
\def\qed{\ifhmode\unskip\nobreak\fi\ifmmode\ifinner\else\hskip5pt
\fi\fi\hbox{\hskip5 pt \vrule width4 pt height6 pt depth1.5 pt \hskip1pt }}

\begin{document}
\title[]{Classification theorems for biharmonic real hypersurfaces in a complex projective space}
\author[]{Toru Sasahara}
\address{
Center for Liberal Arts and Sciences, 
Hachinohe Institute of Technology, 
Hachinohe, Aomori 031-8501, Japan}
\email{sasahara@hi-tech.ac.jp}


\begin{abstract} 
First, we classify proper biharmonic Hopf real hypersurfaces in $\mathbb{C}P^2$.
 Next, we classify proper biharmonic real hypersurfaces with  two distinct principal curvatures 
 in $\mathbb{C}P^n$, where $n\geq 2$.
 Finally, we prove that biharmonic ruled real hypersurfaces in  $\mathbb{C}P^n$
 are minimal, where $n\geq 2$.
\end{abstract}



\maketitle

 \section{Introduction}

Let $f:M \rightarrow N$ be a smooth map between  two Riemannian manifolds.
The bienergy functional $E_2(f)$ of $f$ over compact domain $\Omega\subset M$ is defined by 
\be
E_2(f)=\int_{\Omega}|\tau(f)|^2dv,\nonumber
\ee
where
$
\tau(f):=
\sum_{i=1}^{m}\{\nabla^{f}_{e_i}df(e_i)
-df(\nabla_{e_i}e_i)\}
$ is the tension field of $f$.
Here $\nabla^f$, $\nabla$ and $\{e_i\}_{i=1}^{m}$
 denote the induced connection, the connection of $M$ and
a local orthonormal basis of $M$, respectively. 
%
If $f$ is a critical point of $E_2$
under compactly supported variations, then $f$ is called a
{\it biharmonic map}. 
Jiang \cite{ji2} proved that $f$ is biharmonic if and only if its bitension field defined by
\be
\tau_2(f):=\sum_{i=1}^{m}\left\{(
\nabla^{f}_{e_i}\nabla^{f}_{e_i}-\nabla^{f}_
{\nabla_{e_i}e_i})\tau(f)+R^{N}(\tau(f),df(e_i))df(e_i)\right\}\nonumber
\ee 
vanishes identically, 
where $R^{N}$ is 
the curvature tensor  of $N$,
which is given by $$R^{N}(X, Y)Z=[\nabla^N_X, \nabla^N_Y]Z-\nabla^N_{[X, Y]}Z$$ for the Levi-Civita connection $\nabla^N$ of $N$.

A submanifold is called a {\it biharmonic submanifold} if the isometric immersion that defines the submanifold
is a biharmonic map. Clearly, any minimal submanifold is  biharmonic. Non-minimal biharmonic submanifolds are said to be {\it proper}.
Considerable advancement
has been made in the study
of proper
biharmonic submanifolds 
in manifolds with special metric properties (e.g., real space forms, complex space forms, Sasakian
space forms, conformally flat spaces, etc.) since the beginning of this century.
This paper deals with proper biharmonic real hypersurfaces in 
the complex projective space  $\mathbb{C}P^n(4)$ of complex dimension $n$ and constant
holomorphic  sectional curvature $4$.

Let $M$ be a  real hypersurface in  $\mathbb{C}P^n(4)$.   We denote by $J$ the almost 
complex structure of $\mathbb{C}P^n(4)$. 
The characteristic vector field on $M$ is defined by $\xi=-JN$ for a  local unit normal vector field $N$.
A real hypersurface $M$ is said to be {\it Hopf} if $\xi$ is a principal curvature vector at every point of $M$.
Let $\mathcal{H}$ be the holomorphic distribution defined by $\mathcal{H}=\bigcup_{p\in M}\{X\in T_pM\ |\ \<X, \xi\>=0\}$.
If $\mathcal{H}$ is integrable and each leaf of its maximal integral manifolds is 
locally congruent to
$\mathbb{C}P^{n-1}(4)$, then $M$ is said to be  {\it ruled}.  
A ruled real hypersurface is not a Hopf hypersurface.

In \cite{IIU1}, Ichiyama, Inoguchi and Urakawa  presented the classification result for proper biharmonic homogeneous
real hypersurfaces in  $\mathbb{C}P^n(4)$. 
We will show that this result is not entirely correct and also provide a complete classification.
We obtain, moreover, the following classification results for biharmonic
real hypersurfaces
in $\mathbb{C}P^n(4)$. 
\begin{theorem}\label{thm1}
Let $M$ be a Hopf  hypersurface in $\mathbb{C}P^2(4)$.
Then $M$ is proper biharmonic if and only if it is congruent to
an open part of a geodesic sphere of radius
\begin{eqnarray*}
r=\cot^{-1}
\sqrt{
\frac
{4\pm
\sqrt{
13
}
}
{3}}.
\end{eqnarray*}
\end{theorem}

  \begin{theorem}\label{thm2}
Let $M$ be a real hypersurface with two distinct principal curvatures in $\mathbb{C}P^n(4)$, where 
$n\geq 2$.
Then $M$ is proper biharmonic if and only if it is congruent to
an open part of a geodesic sphere of radius
\begin{eqnarray*}
r=\cot^{-1}
\sqrt{
\frac
{n+2\pm
\sqrt{
n^2+2n+5
}
}
{2n-1}}.
\end{eqnarray*}
\end{theorem} 

\begin{theorem}\label{thm3}
Let $M$ be a ruled real hypersurface in $\mathbb{C}P^n(4)$, where  $n\geq 2$.
If $M$ is biharmonic, then it is minimal.
\end{theorem}

\begin{remark}
	{\rm  Balmu\c{s},  Montaldo and Oniciuc \cite{bal} classified 
		proper biharmonic hypersurfaces with two distinct principal curvatures in the unit sphere.}
\end{remark}

\begin{remark}
{\rm  Adachi, Bao and  Maeda \cite{ada} classified minimal ruled real hypersurfaces in non-flat complex space forms.}
\end{remark}
Throughout this paper, all hypersurfaces are assumed to be orientable.
 \section{Preliminaries}
 \subsection{Fundamental equations and facts for real hypersurfaces}
Let $M$ be a   real hypersurface in  $\mathbb{C}P^n(4)$. 
We denote by $\nabla$ and
 $\tilde\nabla$ the Levi-Civita connections on $M$ and $\mathbb{C}P^n(4)$, respectively. Let $N$ denote a  unit normal vector field on $M$. For any vector fields $X$, $Y$ tangent to $M$, the
 Gauss and Weingarten formulas are  given by
\be
 \begin{split}
 \tilde \nabla_XY&= \nabla_XY+\<AX, Y\>N, \label{gawe}\\
 \tilde\nabla_X N&= -AX,
 \end{split}\nonumber
\ee
 respectively, where $A$ is the shape operator.
The mean curvature function $H$ is defined by 
$H={\rm tr}A/(2n-1).$
 If it vanishes identically, then $M$ is said to be  {\it minimal}.

For any  vector field  $X$ tangent to $M$,  we denote the tangential component of $JX$ by $\phi X$.
Then, from $\tilde\nabla J=0$ and  the Gauss and  Weingarten formulas, it follows  that
\be
\nabla_{X}\xi=\phi AX. \label{PA}
\ee

 
 We denote by $R$ the Riemannian curvature tensor of $M$. Then,
 the equations of Gauss  and Codazzi are given by
 \begin{align}
 &R(X, Y)Z=c[\<Y, Z\>X-\<X, Z\>Y+\<\phi Y, Z\>\phi X
 -\<\phi X, Z\>\phi Y  \label{ga}\\
 &\hskip60pt -2\<\phi X, Y\>\phi Z]
   +\<AY, Z\>AX-\<AX, Z\>AY,\nonumber\\
  &(\nabla_XA)Y-(\nabla_YA)X=c[\<X, \xi\>\phi Y-\<Y, \xi\>\phi X-2\<\phi X, Y\>\xi],\label{co}
 \end{align}       
 respectively.
For later use we need the following fundamental results for real hypersurfaces in $\mathbb{C}P^n(4)$.

 \begin{theorem}[\cite{ni}]\label{eigen}
 Let $M$ be a Hopf hypersurface in $\mathbb{C}P^n(4)$ 
 with $A\xi=\delta\xi$.
 Then
 \begin{itemize}
 \item[(i)] $\delta$ is constant{\rm ;}
 \item[(ii)] If $X$ is a  tangent vector of $M$ orthogonal to $\xi$  such that 
  $AX=\lambda_1X$ and $A\phi X=\lambda_2\phi X$, then
 $2\lambda_1\lambda_2=(\lambda_1+\lambda_2)\delta+2$ holds.
 \end{itemize} 
 \end{theorem}
 
 
\begin{theorem}[\cite{ni}]\label{two}
Let $M$ be a real hypersurface in $\mathbb{C}P^n(4)$, where $n\geq 3$.
If M has two  distinct principal curvatures, then $M$ is an open part of a geodesic sphere.
\end{theorem} 
 
 \begin{theorem}[\cite{ivey}]\label{hopf}
 Let $M$ be a Hopf hypersurface in $\mathbb{C}P^2(4)$ with two distinct principal curvatures.
 Then $M$ is an open part of  a geodesic sphere.
 \end{theorem}

 \begin{lemma}[\cite{ni}]\label{lem2}
Let $M$ a  real hypersurface $M$ in $\mathbb{C}P^n(4)$ with $n\geq 2$.
We define differentiable functions $\alpha$, $\beta$ on $M$ by $\alpha=\<A\xi, \xi\>$
and $\beta=\|A\xi-\alpha\xi\|$.
Then, $M$ is ruled if and only if 
the following two conditions hold{\rm :} 




{\rm (1)} the set $M_1=\{p\in M \ | \ \beta(p)\ne 0\}$ is an open dense subset of $M${\rm ;}

{\rm (2)} there is a unit vector field $U$ on $M_1$, which is orthogonal to $\xi$ and satisfies
\be 
A\xi=\alpha\xi+\beta U, \ \ AU=\beta\xi, \ \ AX=0 \label{ruled}
\ee
 for an arbitrary tangent  vector $X$ orthogonal to both $\xi$ and $U$.

\end{lemma}

We have the following characterization for a real hypersurface in $\mathbb{C}P^n(4)$
 to be
biharmonic.
\begin{proposition}[\cite{fet}]
Let $M$ be a real hypersurface in $\mathbb{C}P^n(4)$.
Then $M$ is biharmonic if and only if the following two equations hold{\rm :}
\begin{align}
&\Delta H+(\|A\|^2-2(n+1))H=0, \label{bihar1}\\
&2A(\nabla H)+(2n-1)H\nabla H=0.\label{bihar2}
\end{align}
\end{proposition} 

\subsection{Proper biharmonic homogeneous real hypersurfaces in $\mathbb{C}P^n(4)$}

In \cite[Theorem 14]{IIU1}, Ichiyama, Inoguchi and Urakawa
presented the classification result for proper biharmonic
homogeneous real hypersurfaces in  $\mathbb{C}P^n(4)$.
However, it is incorrect. In this subsection, we  provide a correct classification.

The following is  a
 complete list of homogeneous real hypersurfaces in $\mathbb{C}P^n(4)$ (see \cite{ni}, \cite{takagi}).
\begin{itemize}
\item {\rm (A)} A tube of radius $r$ over  a totally geodesic $\mathbb{C}P^m(4)$ $(0\leq m\leq n-2)$, where
$0<r<\pi/2$.

\item {\rm (B)} A tube of radius $r$ over a complex quadric $Q^{n-1}$, where
$0<r<\pi/4$.

\item {\rm (C)} A tube of radius $r$ over a $\mathbb{C}P^1(4)\times 
\mathbb{C}P^{(n-1)/2}(4)$, $0<r<\pi/4$ and $n$ $(\geq 5)$ is odd.

\item {\rm (D)} A tube of radius $r$ over the Pl\"{u}cker embedding of the Grassmanian ${\rm Gr}_2(\mathbb{C}^5)\subset\mathbb{C}P^9(4)$, 
where $0<r<\pi/4$.

\item {\rm (E)} A tube of radius $r$ over the canonical embedding of $SO(10)/U(5)\subset\mathbb{C}P^{15}(4)$ 
where $0<r<\pi/4$.
\end{itemize}

On page 251 in \cite{IIU1}, the constant $\|A\|^2$ of real hypersurface of type $D$
 should be
replaced by
\be 
\frac{5X^4-4X^3+62X^2-4X+5}{X(X-1)^2},\nonumber\\
\ee
where $X=\cot^2r$. In addition, on page 252, the constant $\|A\|^2$ of real hypersurface of type $E$
should be
replaced by
\be 
\frac{3(3X^4-2X^3+30X^2-2X+3)}{X(X-1)^2}-2.\nonumber
\ee

Accordingly, equation (\ref{bihar1}) implies that 
the conditions 
for  real hypersurfaces of type $D$ and type $E$ to be proper biharmonic  are given by
\bea
&& 5X^4-24X^3+102X^2-24X+5=0, \label{D}\\
&& 9X^4-40X^3+158X^2-40X+9=0, \label{E}
\eea
respectively.
The left-hand sides of (\ref{D}) and (\ref{E}) can be rewritten, respectively,  as
\bea
&& X^2(5X^2-24X+51)+51X^2-24X+5,\nonumber\\
&& X^2(9X^2-40X+79)+79X^2-40X+9,\nonumber
\eea
 which are positive, because polynomials $5X^2-24X+51$, $51X^2-24X+5$, $9X^2-40X+79$ and
  $79X^2-40X+9$ are positive.
 Therefore, equations (\ref{D}) and (\ref{E}) have no real solutions, and hence
 there exists no proper biharmonic real hypersurface of type $D$ or type $E$
 in $\mathbb{C}P^n(4)$.

Consequently, the correct statement of (II) of Theorem 14 in  \cite{IIU1} is as follows.
\begin{theorem}\label{hom}
Let $M$ be a homogeneous real hypersurface in $\mathbb{C}P^n(4)$.
Then $M$ is proper biharmonic if and only if it is congruent to
an open part of a tube over $\mathbb{C}P^m(4)$ $(0\leq m\leq n-2)$
 of radius 
\be
r=\cot^{-1}\sqrt{\frac{n+2\pm\sqrt{(2m-n+1)^2+4(n+1)}}{2n-2m-1}}. 
\ee
\end{theorem}

\begin{remark}\label{rem1}
	{\rm A real hypersurface in $\mathbb{C}P^n(4)$ is homogeneous if and only if
	 it is a Hopf hypersurface with constant principal curvatures (see \cite{kim}).}
\end{remark}


\section{Proof of Theorem 1.1}
Let $M$ be a Hopf  hypersurface in $\mathbb{C}P^2(4)$.
We assume that $H$ is not constant. Then there exists an open subset $\mathcal{U}$ of $M$
on which $\nabla H\ne 0$. By Lemma 2.13 in \cite{ni}, we have $\xi H=0$, which implies that
$\nabla H$ is normal to $\xi$. Thus, on $\mathcal{U}$ we can choose a local orthonormal frame 
$\{e_1, e_2, e_3\}$ such that
$e_1=\nabla H/\|\nabla H\|$, $e_2=\phi e_1$ and $e_3=\xi$. 
Since $M$ is biharmonic, by  (\ref{bihar2}) we have
 \bea
 A(\nabla H)=-(3H/2)(\nabla H),\nonumber
 \eea
from which it follows  that
the shape operator of $M$ takes the form
\bea
A=
    \begin{pmatrix}
       -3H/2&  0 & 0\\
       0 & \lambda &0 \\
       0 & 0 & \delta
   \end{pmatrix}\nonumber
\eea
for some functions $\lambda$ and $\delta$. By the relation ${\rm tr}A=3H$ we have
\bea\label{mean}
\lambda+\delta=(9/2)H.
\eea

On the other hand, by (ii) of Theorem \ref{eigen}, we are led to
\bea 
-6\lambda H=(2\lambda-3H)\delta+4.\nonumber
\eea
Combining this with (\ref{mean}) and taking into account  (i) of Theorem \ref{eigen},
we see
 that $H$ must be  constant on $\mathcal{U}$. This contradicts our assumption. Therefore, $H$ is constant on $M$. The constancy of ${\rm tr}A$ and Theorem \ref{eigen} show that $M$ has constant principal curvatures.
 Taking into account  Remark \ref{rem1} and applying Theorem \ref{hom},
  we conclude.
  \qed

\section{Proof of Theorem 1.2}
Let $M$ be a real hypersurface with two distinct principal curvatures in $\mathbb{C}P^n(4)$.
If $M$ is a 3-dimensional Hopf hypersurface or if $n\geq 3$, then, 
using Theorem 2.4 for $m=0$, and Theorems 2.2 and 2.3, it is easy to see that
Theorem 1.2 holds.

Thus, we need only consider  the case where $M$ is non-Hopf and  $n=2$.
According to \cite{ivey}, 
on such a real hypersurface with two distinct principal curvatures,
  there exists a local orthonormal frame $\{\xi, X, \phi X\}$ such that
the shape operator takes the form
\bea
A=
    \begin{pmatrix}
       \alpha &  \beta & 0\\
       \beta & \gamma &0 \\
       0 & 0 & \mu
   \end{pmatrix}, 
\eea
where all components of $A$ are constant along the distribution spanned by $\{\xi, X\}$, and 
satisfy
\begin{align}
&\mu^2-(\alpha+\gamma)\mu+\alpha\gamma-\beta^2=0, \label{Eq1} \\
&\phi X(\alpha)=\beta(\alpha+\gamma-3\mu), \label{Eq2}\\
&\phi X(\beta)=\beta^2+\gamma^2+\mu(\alpha-2\gamma)+1, \label{Eq3}\\
&\phi X(\gamma)=\frac{(\gamma-\mu)(\gamma^2-\alpha\gamma-1)}{\beta}
+\beta(2\gamma+\mu).\label{Eq4}
\end{align}

Assume that $M$ is proper biharmonic.
Equation (\ref{bihar2}) for  $n=2$ is equivalent to
\be
(\alpha+\gamma+3\mu)\phi X(H)=0.\label{Eq7}
\ee
We suppose that  there exists an open subset
$\mathcal{U}$ of $M$
on which 
\be
\alpha+\gamma+3\mu=0. \label{Eq5}
\ee
Eliminating $\mu$ from (\ref{Eq1}) and (\ref{Eq5}) leads to
\be
4\alpha^2-9\beta^2+17\alpha\gamma+4\gamma^2=0.\label{poly1}
\ee
By differentiating this equation along  $\phi X$ and  using (\ref{Eq2})-(\ref{Eq4}), we obtain
\be 
\begin{split}
&54\beta^4-(49\alpha^2+209\alpha\gamma+52\gamma^2-54)\beta^2\\
&-32\gamma^4-44\alpha\gamma^3+(59\alpha^2+32)\gamma^2+(17\alpha^3+76\alpha)\gamma
+17\alpha^2=0.
\end{split}\label{Eq6}
\ee
Eliminating $\beta$ from (\ref{poly1}) and (\ref{Eq6}) gives
\be
(\alpha+4\gamma)f(\alpha, \gamma)=0,
\ee
where $f(\alpha, \gamma)$ is given by the following polynomial
\be 
f(\alpha, \gamma):=100\gamma^3+300\alpha\gamma^2
+(300\alpha^2-126)\gamma+100\alpha^3-369\alpha.\nonumber
\ee

If $\alpha+4\gamma=0$, then by (\ref{Eq5}) we have $\gamma=\mu$. 
Combining this and   (\ref{Eq1}) implies $\beta=0$, which does not occur since
$M$ is non-Hopf. Thus, we have $f(\alpha, \gamma)=0$.
Differentiating $f(\alpha, \gamma)=0$ along $\phi X$ and using (\ref{Eq2})-(\ref{Eq4}) gives
\bea
&&1000\gamma^4+2600\alpha\gamma^3+(5000\alpha^2+700\alpha
-1113)\gamma^2+(4400\alpha^3+500\alpha^2-2055\alpha)\gamma
\nonumber\\
&&+1000\alpha^4-200\alpha^3-1842\alpha^2-450\alpha+189=0.\nonumber
\eea
Eliminating $\gamma$ from this equation and $f(\alpha, \gamma)=0$, we get a
non-trivial algebraic equation of 
$\alpha$ with constant coefficients. Thus, $\alpha$ must be a constant.
From (\ref{Eq2}) we obtain $\alpha+\gamma-3\mu=0$, which together with (\ref{Eq1}) and (\ref{Eq5})
shows $\beta=0$. This contradicts the condition that $M$ is non-Hopf.
Therefore, $\mathcal{U}$ is empty, and hence it follows from (\ref{Eq7})  that
$H$ must be constant on $M$.

We put $\alpha+\gamma+\mu=d$ $(\ne0)$ for some non-zero constant $d$. Eliminating $\mu$ from
this equation and (\ref{Eq1}) gives
\be 
2\alpha^2-\beta^2-3\alpha d+d^2+5\alpha\gamma-3d\gamma+2\gamma^2=0.\label{Eq8}
\ee
Differentiating (\ref{Eq8}) along $\phi X$  and  using (\ref{Eq2})-(\ref{Eq4}), we have 
\be
\begin{split}
& 2\beta^4-\{18\gamma^2+(35\alpha-22d)\gamma+13\alpha^2-18\alpha d
+6d^2-2\}\beta^2\\
&-8\gamma^4-(6\alpha-10d)\gamma^3+(9\alpha^2
-2\alpha d-3d^2+8)\gamma^2\\
&+\{5\alpha^3-8\alpha^2d+(3d^2+14)\alpha
-10d\}\gamma+5\alpha^2-8\alpha d+3d^2=0.
\end{split}\label{Eq9}
\ee
Eliminating $\beta$ from (\ref{Eq8}) and (\ref{Eq9}), we obtain
$(\alpha-d+2\gamma)g(\alpha, \gamma)=0$, where
$g(\alpha, \gamma)$ is given by the following polynomial
\be 
\begin{split}
g(\alpha, \gamma):=& 18\gamma^3+(54\alpha-33d)\gamma^2+(54\alpha^2-66\alpha d+20d^2-6)\gamma\\
&+18\alpha^3-33\alpha^2 d+(20d^2-9)\alpha-4d^3+5d.
\end{split}\nonumber
\ee
 If there exists an open subset
of $M$
on which 
$\alpha-d+2\gamma=0$, then we have $\gamma=\mu$. However, this cannot occur because of 
$\beta\ne 0$ and
(\ref{Eq1}).  

Hence, we get $g(\alpha, \gamma)=0$. Differentiating this equation along $\phi X$  and  using (\ref{Eq2})-(\ref{Eq4}) gives
\be
\biggl(\sum_{i=0}^3P_i(\alpha)\gamma^i\biggr)\beta^2+\sum_{i=0}^5Q_i(\alpha)\gamma^i=0, \label{poly2}\nonumber
\ee where $P_i(\alpha)$ and $Q_i(\alpha)$
are polynomials in $\alpha$. 
 Eliminating $\beta$ from this equation and (\ref{Eq8}), we get
 $$(\alpha-d+2\gamma)\sum_{i=0}^4R_i(\alpha)\gamma^i=0,$$
where $R_i(\alpha)$ are polynomials in $\alpha$.
We do not list $P_i$, $Q_i$ and $R_i$ explicitly, however, these polynomials can be recovered quickly by a computer algebra program.

Eliminating $\gamma$ from $g(\alpha, \gamma)=0$ and $\sum_{i=0}^4R_i(\alpha)\gamma^i=0$,
we obtain a non-trivial algebraic equation of $\alpha$ with constant coefficients. Therefore $\alpha$ must be constant, and hence
it follows from (\ref{Eq2}) that $\alpha+\gamma-3\mu=0$,
which together with   (\ref{Eq1}) and $\alpha+\gamma+\mu=d$,
 yields that $\beta$ is also constant, and hence the right-hand side of (\ref{Eq3}) must vanish.
 However, substituting $\mu=(\alpha+\gamma)/3$ into 
 the right-hand side of (\ref{Eq3}), we find that it is positive. Thus, equation (\ref{Eq3}) does not hold. 
 Consequently, we conclude that there exist no  proper biharmonic non-Hopf
real hypersurfaces with two distinct principal curvatures in $\mathbb{C}P^2(4)$. 
 \qed

\section{Proof of Theorem 1.3}
Let $M$ be a ruled real hypersurface in $\mathbb{C}P^n(4)$ whose shape operator
is given by (\ref{ruled}).
Then it follows from (\ref{PA}) that
\bea
&\<\nabla_{\xi}\xi,\phi U\>=\beta, \label{eqq1}\\
&\nabla_{\phi U}\xi=0.\label{eqq5}
\eea
For an arbitrary tangent vector $X$ orthogonal to $\xi$, $U$ and $\phi U$,
 we have (see \cite{kim1})
\be
\nabla_XU=-\beta^{-1}\phi X. \label{eqq2}
\ee
 The following relations have also been derived in \cite{kim1}.
\bea
&\alpha\beta\xi+(\beta^2+1)U 
-A\phi\nabla_{\xi}U-(\phi U\alpha)\xi-(\phi U\beta)U
-\beta\nabla_{\phi U}U=0, \label{eq2} \\
&\nabla_{\phi U}\phi U=0, \quad \nabla_{\phi U}U=0, \label{eq5}\\
&\<\nabla_UU, \phi U\>=\beta-\beta^{-1}, \label{eq3}\\
&\phi U\beta=\beta^2+1. \label{eq4}
\eea
Taking the inner product of each side of (\ref{eq2}) and $\xi$, we obtain
\bea
\phi U\alpha=\alpha\beta+\beta\<\nabla_{\xi}U, \phi U\>. \label{eqq3}
\eea
Using (\ref{eqq5}), (\ref{eqq2}) and (\ref{eq5}),
we get
\begin{align}
\<R(\phi U, \xi)U, \phi U\>
&=\<\nabla_{\phi U}\nabla_{\xi}U, \phi U\>+\<\nabla_{\nabla_{\xi}\phi U}U, \phi U\>\nonumber\\
&=\phi U\<\nabla_{\xi}U, \phi U\>+{\rm tr}\<\nabla_{\xi}\phi U, \cdot\>\<\nabla_{\cdot}U, \phi U\>\nonumber\\
&=\phi U\<\nabla_{\xi}U, \phi U\>+\<\nabla_{\xi}\phi U, U\>\<\nabla_UU, \phi U\>
+\<\nabla_{\xi}\phi U, \xi\>\<\nabla_{\xi}U, \phi U\>.  \label{eqq4}
\end{align}

Assume that $M$ is biharmonic. 
By (\ref{ruled}) we have  $H=\alpha/(2n-1)$ and 
\be
A(\nabla H)=\frac{1}{2n-1}\{(\alpha\xi\alpha+\beta U\alpha)\xi+(\beta\xi\alpha)U\}.\nonumber
\ee
Thus,  taking the inner product of each side of (\ref{bihar2}) and $\phi U$, we obtain
\bea
\alpha\phi U\alpha=0,\label{eq1}\nonumber
\eea
which implies that $\alpha=0$  or $\phi U\alpha=0$. 
If  $\phi U\alpha=0$, then it follows from (\ref{eqq3})  and $\beta\ne 0$ that
\bea
\alpha+\<\nabla_{\xi}U, \phi U\>=0.\label{eq7}
\eea
We differentiate (\ref{eq7}) along $\phi U$. Then, since $\phi U\alpha=0$, we have
\bea
\phi U\<\nabla_{\xi}U, \phi U\>=0. \label{eq8}
\eea
Using (\ref{eqq1}), (\ref{eq3}),  (\ref{eq7}) and (\ref{eq8}), we compute
the right-hand side of (\ref{eqq4}). Then we obtain
\be 
\<R(\phi U, \xi)U, \phi U\>=\alpha(2\beta-\beta^{-1}).\nonumber
\ee 
Combining this 
with the equation (\ref{ga}) of Gauss for $X=\phi U$, $Y=\xi$ and $Z=U$,
we deduce that $\alpha=0$ or $\beta^2=1/2$.
However, by (\ref{eq4}) the latter case cannot occur. 
Consequently, we conclude that  $M$ is minimal. \qed


 \end{document}